\documentclass{amsart}

\usepackage{amsmath, amsfonts, mathrsfs}
\usepackage{amssymb,latexsym}
\usepackage{enumerate}
\usepackage{bbm}

\usepackage[top=1.7cm, bottom=1.7cm, left=3.1cm, right=3.1cm, includefoot, heightrounded]{geometry}

\newtheorem{theorem}{Theorem}[section]

\newtheorem{lemma}[theorem]{Lemma}

\theoremstyle{definition}
\newtheorem{definition}[theorem]{Definition}

\numberwithin{equation}{section}

\begin{document}

\title
{A New Class of Locally Convex Spaces with application}

\author{ Sokol Bush Kaliaj }

\address{
Mathematics Department, 
Science Natural Faculty, 
University of Elbasan,
Elbasan, 
Albania.
}

\email{kaliajsokol@gmail.com}

\thanks{}



\begin{abstract} 
In this paper we define the Radon-Nikodym class ($\mathcal{RN}$ class) of locally 
convex topological vector spaces. 
The $\mathcal{RN}$ class is characterized in terms of the Radon-Nikodym theorem 
for vector measures using integrable by seminorm derivatives.  
It is shown that the  $\mathcal{RNP}$ class of all complete Hausdorff locally convex spaces 
possessing the Radon-Nikodym property  
is properly contained in  the $\mathcal{RN}$ class.  
As an application  we present a Radon-Nikodym theorem for multimeasures with respect to 
the Pettis integrable multifunctions in terms of the $\mathcal{RN}$ class. 
\end{abstract}

\maketitle

\section{Introduction and Preliminaries}

Throughout this paper 
$X$ is a Hausdorff  locally convex space with the topology $\tau$ and   
$\mathscr{P}$ the family of all $\tau$-continuous seminorms. 
$(\Omega, \Sigma, \mu)$ is a complete finite measure space; 
$\Sigma^{+} = \{ E \in \Sigma: \mu(E) >0 \}$ and  
$\Sigma^{+}_{E} = \{ A \in \Sigma^{+}: A \subset E \}$, $E \in \Sigma^{+}$.
For any seminorm $p \in \mathscr{P}$,   
we denote by $\widetilde{X}_{p}$ the quotient vector space $X/p^{-1}(0)$, 
by $\varphi_{p}:X \to \widetilde{X}_{p}$
the canonical quotient map,  
by $(\widetilde{X}_{p},\widetilde{p})$ the quotient normed space  
and by $(\overline{X}_{p},\overline{p})$ the completion of 
$(\widetilde{X}_{p},\widetilde{p})$, $\widetilde{p}(\varphi_{p}(x)) = p(x)$. 
For each $p,q \in \mathscr{P}$ such that $p \leq q$ we can define 
the  continuous linear maps: 
\begin{equation*}
\begin{split}
\widetilde{g}_{pq} : \widetilde{X}_{q} \rightarrow \widetilde{X}_{p}, 
\quad 
\widetilde{g}_{pq}(\varphi_{q}(x))=\varphi_{p}(x)
\end{split} 
\end{equation*}
and $\overline{g}_{pq}: \overline{X}_{q} \rightarrow \overline{X}_{p}$ 
as the continuous linear extension of 
$\widetilde{g}_{pq}$ to $\overline{X}_{q}$.  
We denote by 
\begin{equation*}
\begin{split}
\varprojlim \widetilde{g}_{pq}~\widetilde{X}_{q}
\quad\text{and}\quad
\varprojlim~\overline{g}_{pq}~\overline{X}_{q}
\end{split}
\end{equation*}
the projective limits of the family 
$\{(\widetilde{X}_{p}, \widetilde{p}): p \in \mathscr{P}\}$ 
and 
$\{(\overline{X}_{p}, \overline{p}): p \in \mathscr{P}\}$ 
with respect to the linear maps 
$\widetilde{g}_{pq}$ and $\overline{g}_{pq}$   
respectively, c.f. [18, p.52].


\begin{definition}\label{D_pMeasurable}
A function 
$f : \Omega \to X$ is said to be \textit{measurable} if 
\begin{equation*}
\begin{split}
(\forall G \in \tau)[f^{-1}(G) \in \Sigma]. 
\end{split}
\end{equation*}
If $f$ has only a finite set 
$x_{1}, \dotsc, x_{n}$ of values 
then $f$ is said to be a \textit{simple function}; in this case,  
$f$  is measurable if and only if 
\begin{equation*}
\begin{split}
f^{-1}(x_{i}) = \{t \in \Omega : f(t) = x_{i}\} \in \Sigma, 
\quad
i = 1, \dotsc,n. 
\end{split}
\end{equation*}
\end{definition}


\begin{definition}\label{D_Blondia_M}
Let $f: \Omega \to X$ be a function. 
We say that $f$ is 
\begin{itemize}
\item
\textit{$p$-strongly measurable}($p \in \mathscr{P}$) 
if  there exist 
a sequence of measurable simple functions $(f_{n}^{p} : \Omega \to X)$ 
and a measurable set  
$Z_{p} \in \Sigma$ with $\mu(Z_{p})=0$ 
such that  
\begin{equation*}
\begin{split}
\lim_{n \to \infty} p(f_{n}^{p}(t) -f(t)) = 0
\quad\text{for all }t \in \Omega \setminus Z_{p},
\end{split}
\end{equation*}
\item
\textit{measurable by seminorm} 
if for each $p \in \mathscr{P}$ 
the function $f$ is $p$-strongly measurable, 
\item
\textit{strongly measurable} 
if there exist a sequence of measurable simple functions 
$(f_{n}: \Omega \to X)$ and a measurable set $Z \in \Sigma$ with $\mu(Z)=0$ such that
\begin{equation*}
\begin{split} 
\lim_{n \to \infty} f_{n}(t) = f(t) 
\quad\text{for all }t \in \Omega \setminus Z. 
\end{split}
\end{equation*} 
\end{itemize}
It is clear that if $f: \Omega \to X$ is a measurable simple function, 
then it is strongly measurable. 
\end{definition}


\begin{definition}\label{D_Pairwise_theSame.0}
Let $(f_{s})_{s \in \mathscr{P}}$ be a net of functions
$f_{s}: \Omega \to X$. 
We say that $(f_{s})_{s \in \mathscr{P}}$ is a net  
of \textit{pairwise almost the same functions} 
if for each two continuous seminorms 
$p, q \in \mathscr{P}$, $p \leq q$,  
there exists a measurable set $Z_{pq} \in \Sigma$ with $\mu(Z_{pq})=0$ 
such that  
\begin{equation*}
\begin{split}
p(f_{p}(t)-f_{q}(t))= 0
\quad
\text{for all }t \in \Omega \setminus Z_{pq}.
\end{split}
\end{equation*}  
We say that a Hausdorff locally convex space $X$ 
is called an \textit{$S$-locally convex space}, 
if for each 
net $(f_{s})_{s \in \mathscr{P}}$ of pairwise almost the same functions, 
there exists a function 
$f: \Omega \to X$ 
such that for each $s \in \mathscr{P}$ there is measurable 
$Z_{s} \in \Sigma$ with $\mu(Z_{s}) = 0$ such that 
\begin{equation*}
\begin{split}
s(f(t) - f_{s}(t))= 0,  
\quad
\text{for all }t \in \Omega \setminus Z_{s}
\end{split}
\end{equation*}
We denote by $\mathcal{S}$ the class of all complete $S$-locally convex spaces. 
\end{definition}


\begin{definition}\label{D_B_Integral}
Let $f: \Omega \to X$ be a function and let $p \in \mathscr{P}$. 
We say that $f$ is called  
\textit{$p$-Bochner integrable}  
if there exist 
a sequence of measurable simple functions 
$(f_{n}^{p} : \Omega \to X)$ 
and a measurable set 
$Z_{p} \in \Sigma$ with $\mu(Z_{p})=0$ 
such that
\begin{itemize}
\item[(i)]
$p(f_{n}^{p}(t) -f(t)) \to 0$ for every $t \in \Omega \setminus Z_{p}$,
\item[(ii)]
$p(f_{n}^{p}(\cdot) -f(\cdot)) \in L_{1}(\mu)$ and
\begin{equation*}
\begin{split}
\lim_{n \to \infty} 
\int_{\Omega} p(f_{n}^{p}(t) -f(t)) d \mu = 0,
\end{split}
\end{equation*}
where $L_{1}(\mu)$ is the space of  the Lebesgue integrable functions from $\Omega$ to $\mathbb{R}$, 
\item[(iii)]
for each $E \in \Sigma$ there exists a vector $x_{E}^{p}(f) \in X$ such that 
\begin{equation*}
\begin{split}
\lim_{n \to \infty} 
p \left ( x_{E}^{p}(f) - \int_{E} f_{n}^{p}(t) d \mu \right )= 0.
\end{split}
\end{equation*}
\end{itemize}
It is clear that if $f$ is $p$-Bochner integrable, then 
$\varphi_{p} \circ f$ is Bochner integrable and 
\begin{equation*}
\begin{split}
(B) \int_{E} (\varphi_{p} \circ f)(t) d \mu = \varphi_{p}(x_{E}^{p}(f)). 
\end{split}
\end{equation*}
\end{definition}


The following definition is given in [1, Definition 2.4].

\begin{definition}\label{D_Blondia_I}
A function $f: \Omega \to X$ is said to be 
\textit{integrable by seminorm}  
if for each $p \in \mathscr{P}$ there exist 
a sequence of measurable simple functions 
$(f_{n}^{p} : \Omega \to X)$ 
and a measurable set 
$Z_{p} \in \Sigma$ with $\mu(Z_{p})=0$ 
such that  
\begin{itemize}
\item[(i)]
$p(f_{n}^{p}(t) -f(t)) \to 0$ for every $t \in \Omega \setminus Z_{p}$,
\item[(ii)]
$p(f_{n}^{p}(\cdot) -f(\cdot)) \in L_{1}(\mu)$ and
\begin{equation*}
\begin{split}
\lim_{n \to \infty} 
\int_{\Omega} p(f_{n}^{p}(t) -f(t)) d \mu = 0,
\end{split}
\end{equation*}
\item[(iii)]
for each $E \in \Sigma$ there exists a vector $x_{E} \in X$ such that 
\begin{equation*}
\begin{split}
\lim_{n \to \infty} 
p \left ( x_{E} - \int_{E} f_{n}^{p}(t) d \mu \right )= 0.
\end{split}
\end{equation*} 
We call $x_{E}$ the 
\textit{integral by seminorm} of $f$ over $E$ and set $\int_{E}f(t) d \mu = x_{E}$. 
\end{itemize}
It is clear that if $f$ is integrable by seminorm, then $f$ is 
$p$-Bochner integrable for every $p \in \mathscr{P}$. 
\end{definition}

Let $m: \Sigma \to X$ be a vector measure. 
Given $E \in \Sigma$, 
we denote by $\Pi(E)$ the family of all finite partitions of $E$ 
into elements of $\Sigma$.  
The vector measure $m$ is said to be of \textit{bounded variation} if 
$|m|_{p}(\Omega) < +\infty$ for every $p \in \mathscr{P}$, where 
\begin{equation*}
\begin{split}
|m|_{p}(E) = \sup  
\left \{
\sum_{A \in \pi} p(m(A)) : \pi \in \Pi(E) 
\right \};
\end{split}
\end{equation*}
$|m|_{p}$ is an extend real-valued measure. 
Also $m$ is said to be $\mu$-continuous (shortly  $m \ll \mu$) 
if $\mu(Z) =0$, whenever $Z \in \Sigma$ satisfies $\mu(Z) = 0$. 
It is clear that $m \ll \mu$ if and only if 
for every $p \in \mathscr{P}$ we have $\varphi_{p} \circ m \ll \mu$.  
The set
\begin{equation*}
\begin{split}
A_{E}(m) = 
\left \{
\frac{m(A)}{\mu(A)} \in X : A \in \Sigma^{+}_{E} 
\right \}
\end{split}
\end{equation*}
is called the \textit{average range} of $m$ on $E \in \Sigma$;  
It is said that $m$ has \textit{locally relatively compact average range}  
if $A_{E}(m)$ is relatively compact subset of $X$  
whenever $E \in \Sigma^{+}$.  
We say that 
$m$ is the \textit{indefinite integral} 
of an integrable by seminorm function $f: \Omega \to X$ with respect to $\mu$,   
if 
\begin{equation*}
\begin{split}
m(E) = \int_{E} f(t) d \mu, 
\quad
\text{for every }E \in \Sigma;
\end{split}
\end{equation*}
$f$ is called the 
\textit{Radon-Nikodym derivative} of $m$ with respect to $\mu$.  
A locally convex space $X$ is said to have the \textit{Radon-Nikodym property}  
if for every complete finite measure space $(\Omega, \Sigma, \mu)$ 
and for every vector measure $m: \Sigma \to X$ of bounded variation and $m \ll \mu$, 
the vector measure  
$m$ is the \textit{indefinite integral} 
of an integrable by seminorm function $f: \Omega \to X$ 
with respect to $\mu$.  
For the Radon-Nikodym property in locally convex spaces we refer to 
[2], [3], [4], [9] and [11].

The paper is organized as follows. 
In Section 2 we define $B$-locally convex spaces and the $\mathcal{RN}$ class 
of all complete  $B$-locally convex spaces.
In the section 3, we first present the Radon-Nikodym theorem in $B$-locally convex spaces, 
Theorem \ref{TT_Radon_Nikodym}.  
Then a full characterization of $B$-locally convex spaces in terms of the Radon-Nikodym theorem 
is presented by Theorem \ref{T_CharacterizationS}. 
This theorem yields that   
the $\mathcal{RNP}$ class is contained in the $\mathcal{RN}$ class, 
and since there is a Banach space which does not have the Radon-Nikodym property 
it follows that 
$\mathcal{RNP} \neq \mathcal{RN}$. 
In the last section, we present 
a Radon-Nikodym theorem for multimeasures with respect to the Pettis integrable multifunctions 
in a complete $B$-locally convex spaces, 
Theorem \ref{T_Radon_Nikodym_Multimeasure}.

\section{The definition of the $\mathcal{RN}$ class}

To define the $\mathcal{RN}$ class,  
we need to prove two Lemmas \ref{L_Measurable.0} and  \ref{L_B_Integral.0}. 
\begin{lemma}\label{L_Measurable.0}
If  $f: \Omega \to X$ is $p$-strongly measurable and 
$g: \Omega \to X$ is $q$-strongly measurable and $p \leq q$, then 
$g$ and $f-g$ are $p$-strongly measurable 
and the real valued function $p(f-g)$ is measurable. 
\end{lemma}
\begin{proof}
By hypothesis and Definition \ref{D_Blondia_M} there exist
sequences of measurable simple functions 
$(f_{n}^{p})$ and $(g_{n}^{q})$ 
such that 
\begin{equation*}
\begin{split}
\lim_{n \to \infty} p(f_{n}^{p}(t) - f(t)) = 0
\quad
\lim_{n \to \infty} q(g_{n}^{q}(t) - g(t)) = 0
\quad
\mu-\text{almost everywhere in }\Omega, 
\end{split}
\end{equation*} 
and since $p(g_{n}^{q}(t) - g(t)) \leq q(g_{n}^{q}(t) - g(t))$ 
it follows that $g$ is $p$-strongly measurable. 
Note that
\begin{equation*}
\begin{split}
\lim_{n \to \infty} p[(f_{n}^{p}(t) - g_{n}^{q}(t)) - (f(t)-g(t))] = 0
\quad
\mu-\text{almost everywhere in }\Omega. 
\end{split}
\end{equation*}
This means that $f-g$ is $p$-strongly measurable, 
Hence, $\varphi_{p} \circ (f-g): \Omega \to \overline{X}_{p}$ is strongly measurable. 
Consequently, the real valued function 
$\widetilde{p}[\varphi_{p} \circ (f-g)]$ is measurable, 
and since
\begin{equation*}
\begin{split}
\widetilde{p}[\varphi_{p}(f(t)-g(t))]
=
p(f(t)-g(t))
\end{split}
\end{equation*} 
it follows that the real valued function  
$p(f-g)$ is measurable and this ends the proof.
\end{proof}

\begin{lemma}\label{L_B_Integral.0}
If  $f: \Omega \to X$ is $p$-Bochner integrable and 
$g: \Omega \to X$ is $q$-Bochner integrable and $p \leq q$, then 
$g$ is $p$-Bochner integrable 
and $p(f-g) \in L_{1}(\mu)$. 
\end{lemma}
\begin{proof}
By hypothesis there exist
sequences of measurable simple functions 
$(f_{n}^{p})$, $(g_{n}^{q})$ and the vectors 
$x_{E}^{p}(f), x_{E}^{q}(g)$ satisfying Definition \ref{D_B_Integral}. 

Since each $g_{n}^{q}$ is $p$-strongly measurable and 
$g$ is $q$-strongly measurable we obtain by Lemma \ref{L_Measurable.0} that 
each $p(g_{n}^{q}(\cdot) - g(\cdot))$ is measurable, 
and since 
$p(g_{n}^{q}(\cdot) - g(\cdot)) \leq q(g_{n}^{q}(\cdot) - g(\cdot))$ and 
$q(g_{n}^{q}(\cdot) - g(\cdot)) \in L_{1}(\mu)$ 
it follows that 
$p(g_{n}^{q}(\cdot) - g(\cdot)) \in L_{1}(\mu)$. 
Hence $g$ is $p$-Bochner integrable and
\begin{equation*}
\begin{split}
(B)\int_{E} (\varphi_{p} \circ g)(t) d \mu = \varphi_{p}(x_{E}^{q}(g)),
\quad
\text{for every }E \in \Sigma.
\end{split}
\end{equation*} 
We write
$h_{n}^{p}=f_{n}^{p} - g_{n}^{q}$, $h = f-g$ and note that 
\begin{itemize}
\item[(i)]
$p(h_{n}^{p}(t) -h(t)) \to 0$ at almost all $t \in \Omega$,
\item[(ii)]
$p(h_{n}^{p}(\cdot) -h(\cdot)) \in L_{1}(\mu)$ 
and
\begin{equation*}
\begin{split}
\lim_{n \to \infty} 
\int_{\Omega} p(h_{n}^{p}(t) -h(t)) d \mu = 0,
\end{split}
\end{equation*}
since the real valued function $p(h_{n}^{p}(\cdot) -h(\cdot))$ is measurable and
\begin{equation*}
\begin{split}
\int_{\Omega} p(h_{n}^{p}(t) -h(t)) d \mu 
\leq& 
\int_{\Omega} p(f_{n}^{p}(t) -f(t)) d \mu 
+
\int_{\Omega} p(g_{n}^{q}(t) -g(t)) d \mu\\
\leq& 
\int_{\Omega} p(f_{n}^{p}(t) -f(t)) d \mu 
+
\int_{\Omega} q(g_{n}^{q}(t) -g(t)) d \mu
\end{split}
\end{equation*}
\item[(iii)]
for each $E \in \Sigma$ the vector 
$x_{E}^{p}(h) = x_{E}^{p}(f) + x_{E}^{q}(g)$ is such that 
\begin{equation*}
\begin{split}
\lim_{n \to \infty} 
p \left ( x_{E}^{p}(h) - \int_{E} h_{n}^{p}(t) d \mu \right )= 0,
\end{split}
\end{equation*}  
since
\begin{equation*}
\begin{split} 
p \left ( x_{E}^{p}(h) - \int_{E} h_{n}^{p}(t) d \mu \right ) 
\leq&
p \left ( x_{E}^{p}(f) - \int_{E} f_{n}^{p}(t) d \mu \right ) + 
p \left ( x_{E}^{q}(g) - \int_{E} g_{n}^{q}(t) d \mu \right ) \\
\leq&
p \left ( x_{E}^{p}(f) - \int_{E} f_{n}^{p}(t) d \mu \right ) + 
q \left ( x_{E}^{q}(g) - \int_{E} g_{n}^{q}(t) d \mu \right ).
\end{split}
\end{equation*}
\end{itemize}
This means that $h=f-g$ is $p$-Bochner integrable and
\begin{equation*}
\begin{split}
(B)\int_{E} \varphi_{p}(h(t)) d \mu 
= 
\varphi_{p}(x_{E}^{p}(h))
\quad
\text{for every }E \in \Sigma.
\end{split}
\end{equation*}
Hence, the function
$
p(f(\cdot)-g(\cdot))
=
\widetilde{p}[\varphi_{p}(f(\cdot)-g(\cdot))]
$
is Lebesgue integrable and the proof is finished.
\end{proof}


We are now ready to define the $\mathcal{RN}$ class. 
\begin{definition}\label{D_Pairwise_theSame}
Let $(f_{s})_{s \in \mathscr{P}}$ be a net of 
functions from $\Omega$ to $X$. 
We say that $(f_{s})_{s \in \mathscr{P}}$ is a net  
of \textit{pairwise almost the same Bochner integrable functions} 
if 
$(f_{s})_{s \in \mathscr{P}}$ is a net of 
pairwise almost the same functions   
and each $f_{s}$ is $s$-Bochner integrable. 
We say that a Hausdorff locally convex space $X$ 
is called a \textit{$B$-locally convex space}, 
if for each 
net $(f_{s})_{s \in \mathscr{P}}$ of pairwise almost the same Bochner integrable functions, 
there exists a function 
$f: \Omega \to X$ 
such that 
$s(f(\cdot) - f_{s}(\cdot)) \in L_{1}(\mu)$ and
\begin{equation*}
\begin{split}
\int_{\Omega} s(f(t) - f_{s}(t)) d \mu = 0,  
\quad
\text{for every }s \in \mathscr{P}.
\end{split}
\end{equation*}
We denote by $\mathcal{RN}$ the class of all complete $B$-locally convex spaces. 
It is clear that 
$\mathcal{RNP} \subset \mathcal{RN} \subset \mathcal{S}$. 
\end{definition}

\section{A full characterization of the $\mathcal{RN}$ class}

We first present the Radon Nikodym theorem in a complete $B$-locally convex spaces 
in terms of the average range of a vector measure, 
Theorem \ref{TT_Radon_Nikodym}. 
Then a full characterization of $B$-locally convex spaces in terms of 
the Radon-Nikodym theorem 
is presented in Theorem \ref{T_CharacterizationS}. 
The following auxiliary lemma follows immediately
from [18, II.5.4, p.53] and [18, Exercise 9, p.70].

\begin{lemma}\label{L_SCH} 
Let $(X, \tau)$ be a complete Hausdorff locally convex space and  
let $\mathscr{P}$ be the family of all $\tau$-continuous seminorms.
Then,  
\begin{equation}\label{eq_SCH.1}
L= \varprojlim \widetilde{g}_{pq}~\widetilde{X}_{q}  
=
\varprojlim~\overline{g}_{pq}~\overline{X}_{q} \subset 
\widetilde{X}_{\mathscr{P}} = \prod_{p \in \mathscr{P}} \widetilde{X}_{p} 
\subset \overline{X}_{\mathscr{P}} = \prod_{p \in \mathscr{P}} \overline{X}_{p},
\end{equation}
and the function  
\begin{equation}\label{eq_SCH.2}
\varphi : X \to L, 
\quad
\varphi(x)=(\varphi_{p}(x))_{p \in \mathscr{P}}
\end{equation}
is an isomorphism of $(X,\tau)$ onto $(L, \tau_{L})$, 
where $\tau_{L}$ is the induced topology in $L$ 
by the product topology in $\overline{X}_{\mathscr{P}}$ 
(or by the product topology in $\widetilde{X}_{\mathscr{P}}$). 
\end{lemma}


Lemma \ref{L_SCH} is useful to prove the next auxiliary lemma. 

\begin{lemma}\label{L_Radon_Nikodym}
Let $X$ be a complete Hausdorff locally convex space and 
let $f:\Omega \to X$ be a function. 
Then the following statements are equivalent:
\begin{itemize}
\item[(i)]
$f$ is integrable by seminorm,
\item[(ii)]
$f$ is $p$-Bochner integrable for every $p \in \mathscr{P}$,
\item[(iii)]
$\varphi_{p} \circ f$ is Bochner integrable for every $p \in \mathscr{P}$.
\end{itemize}
In this case, we have
\begin{equation}\label{L_Radon_Nikodym1.0}
\begin{split}
\varphi_{p} \left ( \int_{E} f(t) d \mu  \right ) 
= 
(B)\int_{E} (\varphi_{p} \circ f)(t) d \mu, 
\quad
E \in \Sigma,~p\in \mathscr{P}.  
\end{split}
\end{equation}
\end{lemma}
\begin{proof}
According to Definitions \ref{D_Blondia_I} and \ref{D_B_Integral} 
it follows that 
$(i) \Rightarrow (ii)$ 
and 
$(ii) \Rightarrow (iii)$.

$(iii) \Rightarrow (i)$ 
Assume that $(iii)$ holds. 
Then for each $p \in \mathscr{P}$ there exist 
a sequence of measurable simple functions 
$(\widetilde{g}_{n}^{p} : \Omega \to \widetilde{X}_{p})$ and a measurable set  
$Z_{p} \in \Sigma$ with $\mu(Z_{p})=0$ 
such that  
\begin{itemize}
\item[(i)]
$\widetilde{p}[
\widetilde{g}_{n}^{p}(t) -\varphi_{p}(f(t))] \to 0$ for every $t \in \Omega \setminus Z_{p}$,
\item[(ii)]
$\widetilde{p}[\widetilde{g}_{n}^{p}(\cdot) -\varphi_{p}(f(\cdot))] \in L_{1}(\mu)$ 
and
\begin{equation*}
\begin{split}
\lim_{n \to \infty} 
\int_{\Omega}
\widetilde{p}[
\widetilde{g}_{n}^{p}(t) -\varphi_{p}(f(t))] d \mu = 0. 
\end{split}
\end{equation*}
\end{itemize}
For each $p \in \mathscr{P}$ we can define 
a  measurable simple function $f^{p}_{n}: \Omega \to X$  
by choosing 
$f^{p}_{n}(t) \in \widetilde{g}^{p}_{n}(t)$. 
It is easy to see that sequence $(f^{p}_{n})$ satisfies conditions $(i)$ and $(ii)$ in 
Definition \ref{D_Blondia_I}. 
It remains to check the condition $(iii)$. 
To see this we consider two continuous seminorms 
$p, q \in \mathscr{P}$ such that $p \leq q$.   
It follows by [10, Hille's Theorem II.2.6, p.47]{DIES} that 
\begin{equation*}
\begin{split}
\overline{g}_{pq}
\left ( (B)\int_{E} \varphi_{q}(f(t)) d \mu \right ) 
=&  
(B)\int_{E} (\overline{g}_{pq} \circ \varphi_{q})(f(t)) d \mu \\
=&
(B)\int_{E} (\widetilde{g}_{pq} \circ \varphi_{q})(f(t)) d \mu 
=
(B)\int_{E} \varphi_{p}(f(t)) d \mu.
\end{split}
\end{equation*}
Then,
\begin{equation*}
\begin{split}
\left ( (B)\int_{E} \varphi_{p}(f(t)) d \mu \right )_{p\in \mathscr{P}} \in 
L= \varprojlim \overline{g}_{pq}~\overline{X}_{q}, 
\quad
E \in \Sigma.
\end{split}
\end{equation*}
Consequently,  by Lemma \ref{L_SCH}   
there exists $x_{E} \in X$ such that 
$\varphi_{p}(x_{E}) = (B)\int_{E} \varphi_{p}(f(t)) d \mu$ 
for every $p \in \mathscr{P}$. 
This means that 
$(f^{p}_{n})$ satisfies also the condition $(iii)$ in 
Definition \ref{D_Blondia_I}. 
We now infer that 
$f$ is integrable by seminorm. 
By  equality $\varphi_{p}(x_{E}) = (B)\int_{E} \varphi_{p}(f(t)) d \mu$ 
we obtain \eqref{L_Radon_Nikodym1.0} and the proof is over.
\end{proof}


We now present the Radon-Nikodym theorem in a complete $B$-locally convex space.

\begin{theorem}[Radon-Nikodym theorem]\label{TT_Radon_Nikodym}
Let $X$ be a complete $B$-locally convex space and 
let $m: \Sigma \to X$ be a vector measure. 
Then $m$ is the indefinite integral
of an integrable by seminorm function $f: \Omega \to X$ 
if and only if
\begin{itemize}
\item[(i)]
$m \ll \mu$,  
\item[(ii)]
$m$ is of bounded variation,
\item[(iii)]
$m$ has locally relatively compact average range.
\end{itemize}
\end{theorem}
\begin{proof} 
Assume that 
$f$ is integrable by seminorm and 
\begin{equation*}
\begin{split}
m(E) = \int_{E} f(t) d \mu,
\quad
\text{for every }E \in \Sigma.
\end{split}
\end{equation*}
Then, by Lemma \ref{L_Radon_Nikodym} we obtain 
that $(\varphi_{p} \circ f)$ is Bochner integrable and
\begin{equation*}
\begin{split}
(\varphi_{p} \circ m)(E) = (B)\int_{E} (\varphi_{p} \circ f)(t) d \mu,
\quad
\text{for every }E \in \Sigma.
\end{split}
\end{equation*}
Hence, by [17, Main Theorem] it follows  
that $\varphi_{p} \circ m \ll \mu$, $|m|_{p}(\Omega) < +\infty$ 
and $A_{E}(\varphi_{p} \circ m) = \varphi_{p}\left ( A_{E}(m) \right )$ 
is relatively compact subset of $\widetilde{X}_{p}$ $(E \in \Sigma^{+})$.
By [13, Tychonoff's Theorem], 
the set  
$\prod_{p \in \mathscr{P}} A_{E}(\varphi_{p} \circ m)$ 
is a relatively compact subset 
of $\prod_{p \in \mathscr{P}} \widetilde{X}_{p}$, 
and since
\begin{equation*}
\begin{split}
\varphi \left ( A_{E}(m) \right ) \subset L \cap \prod_{p \in \mathscr{P}} A_{E}(\varphi_{p} \circ m)
\end{split}
\end{equation*}
it follows that $A_{E}(m)$ 
is relatively compact subset of $X$ $(E \in \Sigma^{+})$, 
where $\varphi$ and $L$ are defined in Lemma \ref{L_SCH}

Conversely, assume that the conditions $(i), (ii)$ and $(iii)$ hold. 
Then, for each $p \in \mathscr{P}$, we have  
$\varphi_{p} \circ m \ll \mu$, $|m|_{p}(\Omega) < +\infty$ 
and $\varphi_{p} \left ( A_{E}(m) \right ) = A_{E}(\varphi_{p} \circ m)$ 
is relatively compact subset of $\widetilde{X}_{p}$ $(E \in \Sigma^{+})$. 
Hence, by [17, Main Theorem]
there exists 
a Bochner integrable function $\overline{h}_{p}: \Omega \to \overline{X}_{p}$ 
such that 
\begin{equation*}
\begin{split}
(\varphi_{p} \circ m)(E) = (B)\int_{E} \overline{h}_{p}(t) d \mu,
\quad
\text{for every }E \in \Sigma.
\end{split}
\end{equation*}
Since $\overline{h}_{p}$ is Bochner integrable it is Pettis integrable, 
and since $(\varphi_{p} \circ m)(E) \in \widetilde{X}_{p}$ for all $E \in \Sigma$ 
it follows that $\overline{h}_{p}(t) \in \widetilde{X}_{p}$ at almost all 
$t \in \Omega$, c.f. [7, Lemma 3.7]. 
Hence, there exists a Bochner integrable function $\widetilde{g}_{p} : \Omega \to \widetilde{X}_{p}$ 
such that $\widetilde{g}_{p}(t) = \overline{h}_{p}(t)$   
at almost all $t \in \Omega$.
Further 
there exist 
a sequence of measurable simple functions 
$(\widetilde{g}_{n}^{p} : \Omega \to \widetilde{X}_{p})$
and 
$Z_{p} \in \Sigma$ with $\mu(Z_{p})=0$ 
such that  
\begin{itemize}
\item[(i)]
$\widetilde{p}(
\widetilde{g}_{n}^{p}(t) -\widetilde{g}_{p}(t)) \to 0$ for every $t \in \Omega \setminus Z_{p}$,
\item[(ii)]
$\widetilde{p}(\widetilde{g}_{n}^{p}(\cdot) -\widetilde{g}_{p}(\cdot)) \in L_{1}(\mu)$ 
and
\begin{equation*} 
\begin{split}
\lim_{n \to \infty} 
\int_{\Omega}\widetilde{p}(
\widetilde{g}_{n}^{p}(t) -\widetilde{g}_{p}(t)) d \mu = 0. 
\end{split}
\end{equation*}
\end{itemize}
Thus for every $E \in \Sigma$ we have 
\begin{equation*}
\begin{split}
(\varphi_{p} \circ m)(E) = (B)\int_{E} \overline{h}_{p}(t) d \mu
=(B)\int_{E} \widetilde{g}_{p}(t) d \mu = 
\lim_{n \to \infty} \int_{E} \widetilde{g}^{p}_{n}(t) d \mu.
\end{split}
\end{equation*} 
For each $p \in \mathscr{P}$ we can define 
a function $g_{p}:\Omega \to X$ by choosing $g_{p}(t) \in \widetilde{g}_{p}(t)$, 
and a measurable simple function $g^{p}_{n}: \Omega \to X$  
by choosing 
$g^{p}_{n}(t) \in \widetilde{g}^{p}_{n}(t)$. 
It is clear that each $g_{p}$ is a $p$-Bochner integrable function and 
\begin{equation}\label{eqT_Radon_Nikodym.1}
\begin{split}
(\varphi_{p} \circ m)(E) = (B)\int_{E} (\varphi_{p} \circ g_{p})(t) d \mu
=
(B)\int_{E} \widetilde{g}_{p}(t) d \mu,
\quad
\text{for every }E \in \Sigma.
\end{split}
\end{equation}

We are going to prove that $(g_{s})_{s \in \mathscr{P}}$ is a net of pairwise almost the 
same Bochner integrable functions.  
To see this we consider two continuous seminorms $p, q \in \mathscr{P}$ such that 
$p \leq q$. 
It follows by [10, Hille's Theorem II.2.6, p.47] 
and \eqref{eqT_Radon_Nikodym.1} that 
\begin{equation*}
\begin{split}
\widetilde{g}_{pq} [(\varphi_{q} \circ m)(E)] 
=& 
(B)\int_{E} \widetilde{g}_{pq} \circ (\varphi_{q} \circ g_{q}) d \mu
=
(B)\int_{E} (\varphi_{p} \circ g_{q}) d \mu,
\end{split}
\end{equation*}
and since 
$\widetilde{g}_{pq} [(\varphi_{q} \circ m)(E)] = (\varphi_{p} \circ m)(E)$ 
we obtain 
\begin{equation*}
\begin{split}
(B)\int_{E} (\varphi_{p} \circ g_{q}) d \mu =(B)\int_{E} (\varphi_{p} \circ g_{p}) d \mu, 
\quad\text{for every }E \in \Sigma.
\end{split}
\end{equation*}
Hence, by [10, Corollary II.2.5, p.47] there exists 
$Z_{pq} \in \Sigma$ with $\mu(Z_{pq}) = 0$ such that  
\begin{equation*}
\begin{split}
(\varphi_{p} \circ g_{q})(t) 
=
(\varphi_{p} \circ g_{p})(t),
\text{ for every }t \in \Omega \setminus Z_{pq}, 
\end{split}
\end{equation*}
and consequently
\begin{equation*}
\begin{split}
p(g_{q}(t)-g_{p}(t)) = 
\widetilde{p}((\varphi_{p} \circ g_{q})(t)
- 
(\varphi_{p} \circ g_{p})(t))=0, 
\text{ for every }t \in \Omega \setminus Z_{pq}.
\end{split}
\end{equation*}
This means that 
$(g_{s})_{s \in \mathscr{P}}$ is a net of pairwise almost 
the same Bochner integrable functions.  
Therefore, by Definitions \ref{D_Pairwise_theSame.0} and \ref{D_Pairwise_theSame}  
there exists a function $f: \Omega \to X$ such that 
$s(f(\cdot) - g_{s}(\cdot)) \in L_{1}(\mu)$ and
\begin{equation*}
\begin{split}
s(f(t) - g_{s}(t)) = 0
\quad\text{for all}\quad t \in \Omega \setminus D_{s}, ~s \in \mathscr{P},
\end{split}
\end{equation*}
or equivalently
\begin{equation*}
\begin{split}
(\varphi_{s} \circ f)(t) = (\varphi_{s} \circ g_{s})(t) = \widetilde{g}_{s}(t) 
\quad\text{for all}\quad t \in \Omega \setminus D_{s}, ~s \in \mathscr{P},
\end{split}
\end{equation*} 
where 
$D_{s} \in \Sigma$ with $\mu(D_{s}) = 0$. 
Consequently each function $\varphi_{s} \circ f$ is Bochner integrable and 
\begin{equation*}
\begin{split}
(B)\int_{E} (\varphi_{s} \circ f)(t) d \mu
=
(B)\int_{E} \widetilde{g}_{s}(t) d \mu = 
(\varphi_{s} \circ m)(E),
\quad
\text{for every }E \in \Sigma, ~s \in \mathscr{P}.
\end{split}
\end{equation*} 
Hence by Lemma \ref{L_Radon_Nikodym} it follows that $f$ is integrable by seminorm  
and 
\begin{equation*}
\begin{split}
m(E) = \int_{E} f(t) d \mu
\quad
\text{for every }E \in \Sigma. 
\end{split}
\end{equation*}
Thus $m$ is the indefinite integral of $f$ 
and this ends the proof. 
\end{proof}

We now present a full characterization of $B$-locally convex spaces in terms 
of the Radon-Nikodym Theorem.

\begin{theorem}\label{T_CharacterizationS}
Let $X$ be a complete Hausdorff locally convex space. 
Then the following statements are equivalent:
\begin{itemize}
\item[(i)]
$X$ is a $B$-locally convex space,
\item[(ii)]
for every vector measure $m: \Sigma \to X$ satisfying conditions 
$(i)-(iii)$ in Theorem \ref{TT_Radon_Nikodym} 
it follows that   
$m$ is the indefinite integral
of an integrable by seminorm function $f: \Omega \to X$.
\end{itemize}
\end{theorem}
\begin{proof}
$(ii) \Rightarrow (i)$ 
Assume that $(ii)$ holds and let $(f_{s})$ be a net of almost the same Bochner integrable functions 
and let $x_{E}^{s} \in X$ be the vector satisfying Definition \ref{D_Pairwise_theSame} for 
$E \in \Sigma$ and $s \in \mathscr{P}$. 
Then
\begin{equation*}
\begin{split}
\varphi_{s}(x_{E}^{s}) = (B)\int_{E} (\varphi_{s} \circ f_{s})(t) d \mu,
\quad
\text{for every }E \in \Sigma,~s \in \mathscr{P}, 
\end{split}
\end{equation*}

$\bullet$  
We are going to prove that $(x_{E}^{s})_{s \in \mathscr{P}}$ is a Cauchy net in $X$. 
Given two continuous seminorms $p, q \in \mathscr{P}$ such that $p \leq q$ 
we obtain by Lemma \ref{L_B_Integral.0} that
\begin{equation*}
\begin{split} 
0 =&
\int_{E} p(f_{p}(t) - f_{q}(t)) d \mu 
= 
\int_{E} \widetilde{p}[(\varphi_{p} \circ f_{p})(t) - (\varphi_{p} \circ f_{q})(t)] d \mu \\
\geq&
\widetilde{p}
\left (
(B)\int_{E} (\varphi_{p} \circ f_{p})(t) d \mu - (B)\int_{E} (\varphi_{p} \circ f_{q})(t) d \mu
\right )
\geq 0.
\end{split}
\end{equation*}
Hence, 
\begin{equation*} 
\begin{split}
p(x_{E}^{p} - x_{E}^{q}) = 
\widetilde{p} 
\left (
(B)\int_{E} (\varphi_{p} \circ f_{p})(t) d \mu - (B)\int_{E} (\varphi_{p} \circ f_{q})(t) d \mu
\right )
=0
\end{split}
\end{equation*}
and since
\begin{equation*}
\begin{split}
r \leq p \leq q  \Rightarrow  r(x_{E}^{p} - x_{E}^{q}) \leq p(x_{E}^{p} - x_{E}^{q}) =0 
\Rightarrow 
r(x_{E}^{p} - x_{E}^{q}) = 0
\end{split}
\end{equation*}
it follows that 
$(x_{E}^{s})_{s \in \mathscr{P}}$ is a Cauchy net in the complete space $X$.
Then we can define a function
\begin{equation*}
\begin{split}
m : \Sigma \to X, \quad m(E) = \lim_{s} x_{E}^{s}.
\end{split}
\end{equation*}
Hence,
\begin{equation}\label{eq_Charact_C.0}
\begin{split}
(\varphi_{r} \circ m) (E) 
= 
\lim_{s} \varphi_{r}(x_{E}^{s})
\quad
\text{for every }
r \in \mathscr{P}.
\end{split}
\end{equation}

$\bullet$ 
Fix an arbitrary continuous seminorm $r \in \mathscr{P}$. 
By Lemma \ref{L_B_Integral.0} it follows that 
$(\varphi_{r} \circ f_{s})_{s \geq r}$ is a net in the  space 
$L_{1}(\mu, \overline{X}_{r})$ 
of Bochner integrable functions from $\Omega$ to  $\overline{X}_{r}$.   
Since
\begin{equation*}
\begin{split}
r \leq p \leq q \Rightarrow
\int_{\Omega} 
\widetilde{r}
[(\varphi_{r} \circ f_{p})(t) - (\varphi_{r} \circ f_{q})(t)]
d \mu
=&
\int_{\Omega} 
r(f_{p}(t) - f_{q}(t))d \mu \\
\leq&
\int_{\Omega} 
p(f_{p}(t) - f_{q}(t))d \mu = 0
\end{split}
\end{equation*}
it follows that 
$(\varphi_{r} \circ f_{s})_{s \geq r}$ is a Cauchy net in  
$L_{1}(\mu, \overline{X}_{r})$.  
Hence there exists a Bochner integrable function 
$\overline{h}_{r}: \Omega \to \overline{X}_{r}$ such that 
\begin{equation}\label{eq_Charact_C.1}
\begin{split}
\lim_{s \geq r} \int_{\Omega} \overline{r}[(\varphi_{r} \circ f_{s})(t) - 
\overline{h}_{r}(t)]  d\mu = 0.
\end{split}
\end{equation}
It follows that
\begin{equation*}
\begin{split}
\lim_{s \geq r}  \varphi_{r}(x_{E}^{s})=
\lim_{s \geq r} (B) \int_{E} (\varphi_{r} \circ f_{s})(t) d \mu = 
(B) \int_{E} \overline{h}_{r}(t) d\mu,
\quad\text{for every } E \in \Sigma. 
\end{split}
\end{equation*}
The last equality together with \eqref{eq_Charact_C.0} yields 
\begin{equation}\label{eq_Charact_C.2}
\begin{split}
(\varphi_{r} \circ m)(E) = \lim_{s \geq r} (B)\int_{E} (\varphi_{r} \circ f_{s})(t) d\mu
=(B)\int_{E} \overline{h}_{r}(t) d\mu, 
\quad\text{for every }E \in \Sigma. 
\end{split}
\end{equation}
This means that $\varphi_{r} \circ m$ 
is a vector measure such that $\varphi_{r} \circ m \ll \mu$,  
$|\varphi_{r} \circ m|(\Omega) < +\infty$ and 
it 
has locally relatively compact average range.

$\bullet$ 
Since $r$ was arbitrary, the last result yields that 
$m$ is a vector measure 
satisfying the conditions $(i)-(iii)$ in Theorem \ref{TT_Radon_Nikodym}. 
Then $m$ is the indefinite integral 
of an integrable by seminorm function $f: \Omega \to X$. 
Thus,
\begin{equation*}
\begin{split}
m(E) = \int_{E} f(t) d\mu, 
\quad\text{for every }E \in \Sigma. 
\end{split}
\end{equation*}
The last result together with 
\eqref{eq_Charact_C.2} and Lemma \ref{L_Radon_Nikodym} yields that 
\begin{equation*}
\begin{split}
(\varphi_{r} \circ m)(E) 
= 
\int_{E} (\varphi_{r} \circ f)(t) d\mu 
=
\int_{E} \overline{h}_{r}(t) d\mu
\quad\text{for every }E \in \Sigma,~r \in \mathscr{P}. 
\end{split}
\end{equation*}
Therefore, we obtain by [10, Corollary II.2.5, p.47] that 
$(\varphi_{r} \circ f)(t) = \overline{h}_{r}(t)$ at almost all $t \in \Omega$.
Hence by \eqref{eq_Charact_C.1} we get
\begin{equation*}
\begin{split}
\lim_{s \geq r} \int_{\Omega} r(f_{s}(t) - f(t)) d\mu
=0
\end{split}
\end{equation*}
and since for $s \geq r$ we have
\begin{equation*}
\begin{split}
0 \leq \int_{\Omega} r(f_{r}(t) - f(t) d \mu 
\leq& 
\int_{\Omega} r(f_{r}(t) - f_{s}(t) d \mu  + \int_{\Omega} r(f_{s}(t) - f(t) d \mu  \\
=&
\int_{\Omega} r(f_{s}(t) - f(t) d \mu
\end{split}
\end{equation*}
it follows that $\int_{\Omega} r(f_{r}(t) - f(t) d \mu = 0$. 
This means that $X$ is a $B$-locally convex space.

By virtue of Theorem \ref{TT_Radon_Nikodym} it follows that 
$(i) \Rightarrow (ii)$, 
and this ends the proof.
\end{proof}


\section{An application}

In this section we present a Radon-Nikodym theorem for multimeasures with respect 
to the Pettis integrable multifunctions. 
This theorem is a version of [7, Theorem 3.1] in a complete $B$-locally convex space. 
In paper [7] are presented the Radon-Nikodym theorems for multimeasures with respect 
to the Pettis integrable multifunctions in non-separable locally convex spaces.

We denote by  $2^{X}$ the family of all nonempty subsets of $X$  and by $ck(X)$ the family of all nonempty, 
convex and compact subsets of $X$ is denoted. $X^{*}$ is the topological dual of $X$. 
For any set $C \subset X$ and any $x^{*} \in X^{*}$, we write
\begin{equation*}
\begin{split}
\delta^{*}(x^{*}, C) = \sup \{ x^{*}(x) : x \in C \}.
\end{split}
\end{equation*}

The following definitions are given in [7, Definitions 2.1 and 2.2].

\begin{definition}\label{D_Pettis_Multi}
A multifunction $F: \Omega \to ck(X)$ is called \textit{Pettis integrable} if
\begin{itemize}
\item[(i)]
$\delta^{*}(x^{*}, F)$ is Lebesgue integrable for every $x^{*} \in X^{*}$, 
\item[(ii)]
for each $E \in \Sigma$, there is $M(E) \in ck(X)$ such that
\begin{equation*}
\begin{split}
\delta^{*}
\left ( 
x^{*}, M(E)
\right )
=
\int_{E} \delta^{*}(x^{*}, F) d \mu
\quad
\text{for every }x^{*} \in X^{*},
\end{split}
\end{equation*}
where the function 
$\delta^{*}(x^{*}, F): \Omega \to \mathbb{R}$ is defined by 
$\delta^{*}(x^{*}, F)(t) = \delta^{*}(x^{*}, F(t))$. 
\end{itemize}
We call $M(E)$ the \textit{Pettis integral} of $F$ over $E$ and set $(P)\int_{E} F d \mu = M(E)$. 
\end{definition}
The Pettis integral for multifunctions was first considered by Castaing and Valadier 
[8, Chapter V] 
and has been widely studied in papers 
[5]-[6], [15]-[16] and [12].
The notion of Pettis integrable function $f: \Omega \to X$ 
as can be found in the literature 
(see [10] and [14] for the Banach space case) 
corresponds to Definition \ref{D_Pettis_Multi} 
for $F(t) = \{f(t)\}$ when the integral $(P)\int_{E} F d \mu$ 
is a singleton.
A function $f: \Omega \to X$ is said to be 
a 
\textit{Pettis integrable selector} of a multifunction $F: \Omega \to ck(X)$ if 
$f(t) \in F(t)$ for all $t \in \Omega$ and $f$ is Pettis integrable.

Given a sequence $(A_{n})$ of subsets of $X$,  we write 
$\sum_{n=1}^{+\infty} A_{n}$ to denote the set of all elements of $X$ 
which can be written as the sum of an unconditionally convergent series 
$\sum_{n=1}^{+\infty} x_{n}$, 
where $x_{n} \in A_{n}$ for every $n \in \mathbb{N}$.  

\begin{definition}\label{D_Multimeasure}
A multifunction $M : \Sigma \to 2^{X}$ is called 
a \textit{strong multimeasure} if
\begin{itemize}
\item[(i)]
$M(\emptyset)= \{0\}$, 
\item[(ii)]
for each disjoint sequence $(E_{n})$ in $\Sigma$, 
we have 
\begin{equation*}
\begin{split}
M \left ( \bigcup_{n=1}^{+\infty} E_{n} \right ) 
= 
\sum_{n=1}^{+\infty} M(E_{n}).
\end{split}
\end{equation*}
\end{itemize} 
\end{definition}
We say that the strong multimeasure $M : \Sigma \to 2^{X}$ is 
$\mu$-continuous (shortly $M \ll \mu$) 
if $M(Z) = \{0\}$ whenever $Z \in \Sigma$ satisfies $\mu(Z) =0$. 
The strong multimeasure $M : \Sigma \to 2^{X}$ is called of bounded variation 
if $|M|_{p}(\Omega) < +\infty$ for every $p \in \mathscr{P}$, where 
\begin{equation*}
\begin{split}
|M|_{p}(E) = \sup 
\left \{
\sum_{A \in \pi} ||M(A)||_{p} : \pi \in \Pi(E) 
\right \}
\end{split}
\end{equation*}
and 
\begin{equation*}
\begin{split}
||M(A)||_{p} = \sup \{p(x): x \in A \}.
\end{split}
\end{equation*} 
A selector $m$ of $M$ is a vector-valued function $m : \Sigma \to X$ such that 
$m(E) \in M(E)$ for every $E \in \Sigma$.


Let $\rho :\Sigma \to \Sigma$ be a lifting on $(\Omega, \Sigma, \mu)$, 
c.f. [7].  
We denote by $\Pi_{\rho}$ the family of all finite partitions of $\Omega$ 
into elements of $\rho(\Sigma) \setminus \{\emptyset\}$. 
If we define in $\Pi_{\rho}$ the natural ordering:
\begin{equation*}
\begin{split}
\pi' \succeq \pi'' \Leftrightarrow 
\pi'\text{ is finer than }\pi'',
\end{split}
\end{equation*}
then $\Pi_{\rho}$ is a directed set with respect to "$\succeq$".  
Note that given a net $(x_{\pi})_{\pi \in \Pi_{\rho}}$ with terms in $X$, 
the family 
\begin{equation*}
\begin{split}
\mathcal{F}(x_{\pi}) = \{ \mathcal{F}_{\pi} : \pi \in \Pi_{\rho}\}, 
\quad 
\mathcal{F}_{\pi} = \{ x_{\pi'} \in X : \pi' \geq \pi \} 
\end{split}
\end{equation*}
is a base filter in $X$  
and, we denote by $\mathcal{U}(x_{\pi})$ the ultrafilter on $X$ containing 
$\mathcal{F}(x_{\pi})$. 
It is well known that 
$x_{0} \in X$ is a limit point of the net 
$(x_{\pi})_{\pi \in \Pi_{\rho}}$ 
if and only if 
the corresponding ultrafilter $\mathcal{U}(x_{\pi})$ 
converges to $x_{0}$; in this case we write
\begin{equation*}
\begin{split}
\lim_{\pi} \mathcal{U}(x_{\pi}) = x_{0}.
\end{split}
\end{equation*}  
Given a net $(V_{\pi})_{\pi \in \Pi_{\rho}}$ with terms in $2^{X}$, 
we write 
\begin{equation*}
\begin{split}
\lim_{\pi} \mathcal{U}(V_{\pi}) = 
\left \{
\lim_{\pi} \mathcal{U}(v_{\pi})  : (\forall \pi \in \Pi_{\rho})[v_{\pi} \in V_{\pi}] 
\right \}
\end{split}
\end{equation*}

\begin{definition}\label{D_MbyS_Selector}
Let $F: \Omega \to ck(X)$ be a multifunction 
and let $f:\Omega \to X$ be a function. 
We say that $f$ is called :
\begin{itemize}
\item[(i)]
\textit{$p$-strongly measurable selector} of $F$ if $f$ is 
$p$-strongly measurable and there exists 
a measurable set $Z_{p} \in \Sigma$ with $\mu(Z_{p}) = 0$ such that 
\begin{equation*}
\begin{split}
(\varphi_{p} \circ f)(t) \in \varphi_{p}[F(t)] = 
\{\varphi_{p}(x): x \in F(t) \}
\quad
\text{for every }
t \in \Omega \setminus Z_{p},  
\end{split}
\end{equation*} 
\item[(ii)]
\textit{measurable by seminorm selector} of $F$, if $f$ is 
measurable by seminorm and $p$-strongly measurable selector of $F$ 
for every $p \in \mathscr{P}$. 
\end{itemize}
We say that 
$f$ is an \textit{integrable by seminorm selector} of $F$ if 
$f$ is integrable by seminorm and 
measurable by seminorm selector of $F$.
\end{definition}

We now present a Radon-Nikodym theorem for multimeasures in 
$B$-locally convex spaces.
\begin{theorem}\label{T_Radon_Nikodym_Multimeasure}
Let $X$ be a complete $B$-locally convex space 
and let $M: \Sigma \to ck(X)$ be a strong multimeasure for which there is a set $Q \in ck(X)$ 
such that 
\begin{equation*}
\begin{split}
M(E) \subset \mu(E) Q
\quad
\text{for all }
E \in \Sigma.
\end{split}
\end{equation*}
Then there exists a Pettis integrable multifunction 
$F: \Omega \to ck(X)$ such that 
\begin{itemize}
\item[(i)]
for every countable additive selector $m$ of $M$ there is an integrable 
by seminorm selector $f$ of $F$ such that $m$ is the indefinite integral of $f$,
\item[(ii)]
for each $E \in \Sigma$ the following equalities hold:
\begin{equation*}
\begin{split}
M(E) = (P)\int_{E} F d \mu = 
I(E)
\quad
\text{for every }E \in \Sigma, 
\end{split}
\end{equation*}
where 
\begin{equation*}
\begin{split}
I(E) = 
\left \{
\int_{E} f d \mu : 
f 
\text{ is an integrable by seminorm selector of }F
\right\}.
\end{split}
\end{equation*}
\end{itemize}
\end{theorem}
\begin{proof} 
$(i)$ 
Let $m$ be a countable additive selector of $M$.
Then $m$ is a vector measure of bounded variation, $m \ll \mu$ 
and $m$ has locally relative compact average range, since 
$m(E) \in \mu(E) Q$ for every $E \in \Sigma$. 
Consequently, by Theorem \ref{T_CharacterizationS} 
there is an integrable by seminorm function 
$f: \Omega \to X$ such that 
\begin{equation*}
\begin{split}
m(E) = \int_{E} f(t) d \mu
\quad\text{for all }E \in \Sigma.
\end{split}
\end{equation*}
For each $\pi \in \Pi_{\rho}$ we write 
\begin{equation*}
\begin{split}
M_{\pi} : \Omega \to ck(X), 
\quad 
M_{\pi}(t) = \sum_{E \in \pi} \frac{M(E)}{\mu(E)} \mathbbm{1}_{E}(t), 
\end{split}
\end{equation*}
and
\begin{equation*}
\begin{split}
m_{\pi} : \Omega \to X, 
\quad 
m_{\pi}(t) = \sum_{E \in \pi} \frac{m(E)}{\mu(E)} \mathbbm{1}_{E}(t).
\end{split}
\end{equation*}

By virtue of [7, Theorem 3.1] there exists a Pettis integrable multifunction 
$F: \Omega \to ck(X)$ such that 
\begin{equation*}
\begin{split}
F(t) = \overline{\lim_{\pi} \mathcal{U}(M_{\pi}(t))}
\quad
\text{for all }t \in \Omega 
\end{split}
\end{equation*}
and
\begin{equation}\label{eq_T_Radon_Nikodym_Multimeasure.0}
\begin{split}
M(E) = 
(P)\int_{E} F d \mu
\quad\text{for every }E \in \Sigma
\end{split}
\end{equation}
and the function
\begin{equation}\label{eq_T_Radon_Nikodym_Multimeasure.1}
\begin{split}
h: \Omega \to X, 
\quad
h(t) = \lim_{\pi} \mathcal{U}(m_{\pi}(t))
\end{split}
\end{equation}
is a Pettis integrable selector of $F$ such that
\begin{equation*}
\begin{split}
m(E) = (P)\int_{E} h(t) d \mu \in M(E).
\end{split}
\end{equation*} 
Moreover $M(E) = P(E)$ for all $E \in \Sigma$, 
where 
\begin{equation*}
\begin{split}
P(E) = 
\left \{
(P)\int_{E} g d \mu : 
g 
\text{ is a Pettis integrable selector of }F
\right\}.
\end{split}
\end{equation*}

Lemma \ref{L_Radon_Nikodym} yields that the function $\varphi_{p} \circ f: \Omega \to \widetilde{X}_{p}$ 
is Bochner integrable and 
\begin{equation*}
\begin{split}
(\varphi_{p} \circ m)(E) = (B)\int_{E} (\varphi_{p} \circ f)(t) d \mu
\quad\text{for all }E \in \Sigma.
\end{split}
\end{equation*}
Hence by [7, Lemma 3.7] there exists 
a measurable set $Z_{p} \in \Sigma$ with $\mu(Z_{p}) = 0$ such that 
\begin{equation*} 
\begin{split}
\lim_{\pi} \varphi_{p}(m_{\pi}(t)) = \lim_{\pi} (\varphi_{p} \circ m)_{\pi}(t) = (\varphi_{p} \circ f)(t)
\quad\text{for all }t \in \Omega \setminus Z_{p},
\end{split}
\end{equation*}
since $\varphi_{p} \circ f$ is also Pettis integrable. 
The last result together with 
\eqref{eq_T_Radon_Nikodym_Multimeasure.1} 
yields  
\begin{equation}\label{eq_T_Radon_Nikodym_Multimeasure.2} 
\begin{split}
(\varphi_{p} \circ f)(t) = (\varphi_{p} \circ h)(t)
\quad\text{for all }t \in \Omega \setminus Z_{p},
\end{split}
\end{equation}
and since 
$(\varphi_{p} \circ h)(t) \in \varphi_{p}[F(t)]$ 
for all $t \in \Omega$ 
it follows that 
$(\varphi_{p} \circ f)(t) \in \varphi_{p}[F(t)]$ 
for all $t \in \Omega \setminus Z_{p}$.
This means that $f$ is a measurable by seminorm selector of $F$. 
Thus, the proof of $(i)$ is completed.

Note that by \eqref{eq_T_Radon_Nikodym_Multimeasure.2} it follows that 
$(\varphi_{p} \circ h)$ is Bochner integrable and  
\begin{equation*} 
\begin{split}
(B)\int_{E} (\varphi_{p} \circ f)(t) d \mu 
= 
(B)\int_{E} (\varphi_{p} \circ h)(t) d \mu
\quad\text{for every }E \in \Sigma, ~p \in \mathscr{P}. 
\end{split}
\end{equation*}
Therefore, by Lemma \ref{L_Radon_Nikodym} we obtain
that $h$ is integrable by seminorm and 
\begin{equation*} 
\begin{split}
m(E) = \int_{E} f(t) d \mu = \int_{E} h(t) d \mu =(P)\int_{E} h(t) d \mu
\quad\text{for every }E \in \Sigma. 
\end{split}
\end{equation*}

$(ii)$ 
Let $h$ be a Pettis integrable selector of $F$ and 
\begin{equation*}
\begin{split}
m_{h}(E) = (P)\int_{E} h(t) d \mu,
\qquad
E \in \Sigma.
\end{split}
\end{equation*}
Then by [7, Theorem 3.1] we have
\begin{equation*}
\begin{split}
m_{h}(E) = (P)\int_{E} h(t) d \mu \in P(E) = M(E). 
\end{split}
\end{equation*} 
Hence by $(i)$ it follows that 
$h$ is an integrable by seminorm selector of $F$. 
Consequently, $m_{h}(E) \in I(E)$. 
Thus, $P(E) \subset I(E)$.

We now assume that $f$ is an integrable by seminorm selector of $F$ and 
\begin{equation*}
\begin{split}
m_{f}(E) = \int_{E} f(t) d \mu,
\qquad
E \in \Sigma.
\end{split}
\end{equation*}
It easy to see that the multifunction  
\begin{equation*}
\begin{split}
\widetilde{F}_{p} : \Omega \to ck(\widetilde{X}_{p}), 
\quad
\widetilde{F}_{p}(t) = \{\varphi_{p}(x): x \in F(t)\}
\end{split}
\end{equation*}
is Pettis integrable and 
by \eqref{eq_T_Radon_Nikodym_Multimeasure.0} it follows that 
\begin{equation*}
\begin{split}
\widetilde{M}_{p}(E) 
= 
(P)\int_{E} 
\widetilde{F}_{p} d \mu,
\qquad
p \in \mathscr{P}, 
\end{split}
\end{equation*}
where 
\begin{equation*}
\begin{split}
\widetilde{M}_{p} : \Sigma \to ck(\widetilde{X}_{p}), 
\quad
\widetilde{M}_{p}(E) = \{\varphi_{p}(x) : x \in M(E)\}.
\end{split}
\end{equation*}
By Definition \ref{D_MbyS_Selector} we have 
$(\varphi_{p} \circ f)(t) \in \widetilde{F}_{p}(t)$ at almost all 
$t \in \Omega$ and each $(\varphi_{p} \circ f)$ is Bochner integrable, 
see Lemma \ref{L_Radon_Nikodym}. 
Hence, $(\varphi_{p} \circ f)$ is Pettis integrable and 
\begin{equation*}
\begin{split}
\varphi_{p}(m_{f}(E)) 
= \int_{E} (\varphi_{p} \circ f)(t)d \mu 
= (P)\int_{E} (\varphi_{p} \circ f)(t)d \mu,
\end{split}
\end{equation*}
and since 
$(P)\int_{E} (\varphi_{p} \circ f)(t)d \mu \in (P)\int_{E} 
\widetilde{F}_{p} d \mu$ 
it follows that
\begin{equation*}
\begin{split}
\varphi_{p}(m_{f}(E)) \in \widetilde{M}_{p}(E)
\quad
\text{for every }p \in \mathscr{P}. 
\end{split}
\end{equation*}
Consequently,
\begin{equation*}
\begin{split}
m_{f}(E) \in M(E),
\qquad
E \in \Sigma. 
\end{split}
\end{equation*}
Thus,
$I(E) \subset M(E), E \in \Sigma$. 
By the last result and  
$P(E) = M(E) \subset I(E)$  
we infer that $I(E) = M(E) = P(E)$, and this ends the proof.
\end{proof}

\newpage

$$
\text{REFERENCES}
$$

$[1]$ C.Blondia, \textit{Integration in locally convex spaces}, Simon Stevin, 55(3), (1981).

$[2]$ C.Blondia, \textit{On the Radon-Nikodym theorem for vector valued measures}, Bull. Soc. Math. Belg. (1981). 

$[3]$ C.Blondia, \textit{On the Radon-Nikodym property in locally convex spaces and the completeness}, (1982).

$[4]$ R.D.Bourgin, \textit{Geometric Aspects of Convex Sets with the Radon-Nikodym Property}, 
Lecture Notes in Mathematics, 993 (1983).  

$[5]$ B.Cascales, V.Kadets, and J.Rodr\'{i}guez, 
\textit{The Pettis integral for multi-valued functions via single-valued ones}, 
J. Math. Anal. Appl. 332 (1) (2007).

$[6]$ B.Cascales, V.Kadets, and J.Rodr\'{i}guez,  
\textit{Measurable selectors and set-valued Pettis integral in non-separable Banach spaces}, 
J. Funct. Anal. 256 (3) (2009).

$[7]$ B.Cascales, V.Kadets, J.Rodriguez,
\textit{Radon-Nikodym theorems for multimeasures in non-sepa-rable spaces}, 
Zh. Mat. Fiz. Anal. Geom. 9 (1) (2013).

$[8]$ C.Castaing, M.Valadier, 
\textit{Convex Analysis and Measurable Multifunctions} ,
Lect. Notes Math., 580 (1977).

$[9]$ G.Y.H.Chi, 
\textit{On the Radon-Nikodym theorem and locally convex spaces with the Radon-Nikodym property}, 
Proc. Am. Math. Soc., 62 (2) (1977).

$[10]$ J.Diestel and J.J.Uhl, 
\textit{Vector Measures}, 
Amer. Math. Soc. (1977).

$[11]$ L.Egghe, 
\textit{On the Radon-Nikodym Property, and related topics in locally
convex spaces}, 
Lect. Notes in Math., 645 (1978).

$[12]$ K.El Amri and C.Hess, 
\textit{On the Pettis integral of closed valued multifunctions}, 
Set-Valued Anal. 8(4) (2000).

$[13]$ J.L.Kelly, 
\textit{General Topology}, 
Springer Verlag, (1955).

$[14]$ K.Musial,
\textit{Topics in the theory of Pettis integration},
Rend. Istit. Mat. Univ. Trieste 23 (1991).

$[15]$ K.Musial,
\textit{Pettis integrability of multifunctions with values in arbitrary Banach spaces}, 
J. Convex Analysis 18 (2011).

$[16]$ K.Musial,
\textit{Approximation of Pettis integrable multifunctions with values in arbitrary Banach spaces}, 
J. Convex Analysis 20(3) (2013).

$[17]$ M.A.Rieffel, 
\textit{The Radon-Nikodym theorem for the Bochner integral}, Trans. Amer. Math. Soc, 131 (1968).

$[18]$ H.H.Schaefer and W.P.Wolff, 
\textit{Topological Vector Spaces},
Originally published by Springer-Verlag New York, (1999).

\end{document}